\documentclass[12pt]{article}
\usepackage{amsmath,amsfonts,latexsym,amsthm,amssymb}
\newcommand{\F}{\mathbf F_q}
\newcommand{\K}{\overline{K}_c}

\DeclareMathOperator{\diag}{diag}
\DeclareMathOperator{\dist}{dist}
\begin{document}
\newtheorem{lem}{Lemma}
\newtheorem{teo}{Theorem}
\pagestyle{plain}
\title{Differential Equations for $\F$-Linear Functions, II:
Regular Singularity}
\author{Anatoly N. Kochubei\footnote{Partially supported by
CRDF under Grants UM1-2090 and UM1-2421-KV-02}
\\ \footnotesize Institute of Mathematics,\\
\footnotesize National Academy of Sciences of Ukraine,\\
\footnotesize Tereshchenkivska 3, Kiev, 01601 Ukraine
\\ \footnotesize E-mail: \ ank@sita.kiev.ua}
\date{}
\maketitle
\vspace*{2cm}
Running head:\quad  ``Regular Singularity''
\newpage
\vspace*{8cm}
\begin{abstract}
We study some classes of equations with Carlitz derivatives for
$\F$-linear functions, which
are the natural function field counterparts of linear ordinary
differential equations with a regular singularity. In particular,
an analog of the equation for the power function, the Fuchs and
Euler type equations, and Thakur's hypergeometric equation are
considered. Some properties of the above equations are similar to
the classical case while others are different. For example, a
simple model equation shows a possibility of existence of a
non-trivial continuous locally analytic $\F$-linear solution
which vanishes on an open neighbourhood of the initial point.
\end{abstract}
\vspace{2cm}
\newpage
\section{INTRODUCTION}

\medskip
Let $K$ be the field of formal Laurent series
$t=\sum\limits_{j=N}^\infty \xi_jx^j$ with coefficients $\xi_j$
from the Galois field $\F$, $\xi_N\ne 0$ if $t\ne 0$, $q=p^\upsilon $,
$\upsilon \in \mathbf Z_+$,
where $p$ is a prime number. It is well known that any non-discrete
locally compact field of characteristic $p$ is isomorphic to such $K$.
The absolute value on $K$ is given by $|t|=q^{-N}$, $|0|=0$. The ring of
integers $O=\{ t\in K:\ |t|\le 1\}$ is compact in the topology
corresponding to the metric $\dist (t,s)=|t-s|$.
Let $\K$ be a completion of an algebraic closure of $K$. The
absolute value $|\cdot |$ can be extended in a unique way onto
$\K$.

A function defined on a $\F$-subspace $K_0$ of $K$, with values
in $\K$, is called $\F$-linear
if $f(t_1+t_2)=f(t_1)+f(t_2)$ and $f(\alpha t)=\alpha f(t)$ for any
$t,t_1,t_2\in K$, $\alpha \in \F$.

Many interesting functions studied in analysis over $K$, like analogs
of the exponential, logarithm, Bessel, and hypergeometric functions
(see e.g. \cite{C1,C2,G1,G2,Th1,Th2,SL,K1,K2}) are $\F$-linear
and satisfy some differential equations with polynomial
coefficients, in which the role of a derivative is played by the operator
$$
d=\sqrt[q]{\ }\circ \Delta ,\quad (\Delta u)(t)=u(xt)-xu(t),
$$
where $x$ is a prime element in $K$. The meaning of a polynomial coefficient
in the function field case is not a usual multiplication by a polynomial,
but the action of a polynomial in the operator $\tau$, $\tau u=u^q$.

This paper is a continuation of the article \cite{K3}, in which
a general theory of such equations was initiated. In particular,
we considered regular equations (or systems) of the form
\begin{equation}
du(t)=P(\tau )u(t)+f(t)
\end{equation}
where for each $z\in \left(\K\right) ^m$, $t\in K$,
\begin{equation}
P(\tau )z=\sum \limits _{k=0}^\infty \pi _kz^{q^k},\quad
f(t)=\sum \limits _{j=0}^\infty \varphi _jt^{q^j},
\end{equation}
$\pi _k$ are $m\times m$ matrices with elements from $\K$,
$\varphi _j\in \left(\K\right) ^m$, and the series (2)
have positive radii of convergence. The action of the operator $\tau $
upon a vector or a matrix is defined component-wise, so that
$z^{q^k}=\tau^k(z)=\left( z_1^{q^k},\ldots ,z_m^{q^k}\right)$ for
$z=(z_1,\ldots ,z_m)$. Similarly, if $\pi =(\pi _{ij})$ is a
matrix, we write $\tau^k(\pi )=\left( \pi_{ij}^{q^k}\right)$.
$\tau$ is an automorphism of the ring of matrices over $\K$.

It was shown in \cite{K3} that for any $u_0\in \left(\K\right)
^m$ the equation (1) has a unique local analytic $\F$-linear
solution satisfying the initial condition
\begin{equation}
\lim\limits_{t\to 0}t^{-1}u(t)=u_0.
\end{equation}
Singular higher order scalar equations were also considered (here
the singularity means that the leading coefficient is a non-constant
analytic function of the operator $\tau$), and it was shown that,
in contrast to the classical analytic theory of differential
equations, any formal power series solution has a positive radius
of convergence. Note however that in general a singular equation
need not possess a formal power series solution.

In this paper we study an analog of the class of equations with
regular singularity, the most thoroughly studied class of singular
equations (see \cite{H} or \cite{CL} for the classical
theory of differential equations over $\mathbb C$; the case of
non-Archimedean fields of characteristic zero was studied in
\cite{DGS,RC}).

A typical class of systems with regular singularity at the origin
$\zeta =0$ over $\mathbb C$ consists of systems of the form
$$
\zeta y'(\zeta )=\left( B+\sum\limits_{k=1}^\infty
A_k\zeta^k\right) y(\zeta )
$$
where $B,A_j$ are constant matrices, and the series converges on
a neighbourhood of the origin. Such a system possesses a
fundamental matrix solution of the form $W(\zeta )\zeta^C$ where
$W(\zeta )$ is holomorphic on a neighbourhood of zero, $C$ is a
constant matrix, $\zeta^C=\exp (C\log \zeta )$ is defined by the
obvious power series. Under some additional assumptions regarding
the eigenvalues of the matrix $B$, one can take $C=B$. For
similar results over $\mathbb C_p$ see Sect. III.8 in \cite{DGS}.

In order to investigate such a class of equations in the
framework of $\F$-linear analysis over $K$, one has to go beyond
the class of functions represented by power series. An analog of
the power function need not be holomorphic, and cannot be defined
as above. Fortunately, we have another option here -- instead of
power series expansions we can use the expansions in Carlitz
polynomials on the compact ring of integers $O\subset K$. It is
important to stress that our approach would fail if we consider
equations over $\K$ instead of $K$ (our solutions may take their
values from $\K$, but they are defined over subsets of $K$). In
this sense our techniques are different from the ones developed
for both the characteristic zero cases.

We begin with the simplest model scalar equation
\begin{equation}
\tau du=\lambda u,\quad \lambda \in \K,
\end{equation}
whose solution may be seen as a function field counterpart of the
power function $t\mapsto t^\lambda$. If $|\lambda |<1$, the
equation (4) has a non-trivial continuous $\F$-linear solution
$u$ on the ring of integers $O$. This
solution is analytic on $O$ if and only if $\lambda =[j]$,
$j=0,1,2,\ldots$; we use Carlitz's notation
$$
[j]=x^{q^j}-x.
$$
In this case $u(t)=ct^{q^j}$, $c\in \K$. If $\lambda \ne [j]$ for any
$j$, then the solution of (4) is locally analytic on $O$ if and
only if $\lambda =-x$, and in the latter case $u(t)\equiv 0$ for
$|t|\le q^{-1}$. This paradoxical fact is a good illustration to
the violation of the principle of analytic continuation in the
non-Archimedean case. It is also interesting that the nonlinear
equation $du=\lambda u$ is much simpler than the linear equation
$\tau du=\lambda u$.

The construction of a solution of the equation (4) is generalized
to the case of systems of equations where $\lambda$ is a matrix.
This makes it possible to study the system
\begin{equation}
\tau du-P(\tau )u=0
\end{equation}
where $P(\tau )$ is a matrix-valued analytic function of $\tau$
of the form (2). We construct a matrix-valued solution of (5)
which is written as $W(g(t))$, where $g(t)$ is a solution of the
equation $\tau dg=\pi_0g$, $W(s)=\sum\limits_{k=0}^\infty
w_ks^{q^k}$ has a non-zero radius of convergence. In contrast to
similar results for equations over $\mathbb C$ \cite{CL,H}, here
we have a composition of matrix-functions instead of their
multiplication. As an example, an Euler-type scalar higher-order
equation is considered.

Finally, we study a problem set by Thakur \cite{Th2}. Within his
theory of hypergeometric functions on $K$, Thakur introduced
\cite{Th1,Th2} an analog of the Gauss hypergeometric function
${}_2F_1(a,b;c;t)$ and the corresponding differential equation.
He constructed two families of analytic solutions which coincide
if $c=1$. Classically (over $\mathbb C$) in the latter case there
is another solution with a logarithmic singularity near the
origin, and a natural question is about the kind of singularity
of a non-analytic solution in our function field case. Note that
the function field hypergeometric equation is nonlinear (it is
only $\F$-linear), and the set of solutions is not parametrized
by parameters from $\K$. We prove that in the case $c=1$ a
generic solution of the hypergeometric equation is defined only
for $t\in \F [x]$, and its formal Fourier-Carlitz series cannot
be extended to non-polynomial arguments.

The author is grateful to D. Thakur for his helpful remarks and
suggestions.

\bigskip
\section{A MODEL SCALAR EQUATION}

\medskip
Let us consider the equation (4). We shall look for a continuous
$\F$-linear solution
\begin{equation}
u(t)=\sum\limits_{i=0}^\infty c_if_i(t),\quad t\in O,
\end{equation}
where $\K\ni c_i\to 0$ as $i\to \infty$, $\{f_i(t)\}$ is the
sequence of orthonormal Carlitz polynomials, that is
$f_i(t)=D_i^{-1}e_i(t)$,
$$
D_0=1,\quad D_i=[i][i-1]^q\cdots [1]^{q^{i-1}},\quad e_0(t)=t,
$$
\begin{equation}
e_i(t)=\prod \limits _{\genfrac{}{}{0pt}{}{\omega \in \mathbf F_q[x]}
{\deg \omega <i}}(t-\omega ),\quad i\ge 1
\end{equation}
(see \cite{G1}). The orthonormality means that
$$
\sup\limits_{t\in O}|u(t)|=\sup\limits_i|c_i|.
$$

In fact, $\{ f_i\}$ is a basis of the space of continuous $\F$-linear
functions on $O$ with values in $\K$ -- if $u$ is such a
function, then it can be represented by a convergent series (6)
where the coefficients $c_i\to 0$ are determined uniquely.
Conversely, a series (6) with $c_i\to 0$ defines a continuous $\F$-linear
function.

The polynomials $f_i$ are the characteristic $p$ analogs of the
binomial coefficients forming the Mahler basis of the space of
continuous functions on $\mathbb Z_p$, the elements $D_i$ are the
counterparts of the factorials $i!$.

The rate of decay of the coefficients $c_i$ corresponds to the
smoothness properties of a function $u$ (see \cite{K2,Y}). In
particular, $u(t)$ is locally analytic if and only if
\begin{equation}
\gamma =\liminf\limits_{n\to \infty }\left\{ -q^{-
n}\log_q|c_n|\right\} >0,
\end{equation}
and if (8) holds, then $u(t)$ is analytic on any ball
of the radius $q^{-l}$, $l=\max (0,[-(\log (q-1)+\log \gamma
)/\log q]+1)$ (see \cite{Y}). Note, in order to avoid confusion,
that the formula (8) looks different from the corresponding
formula in Sect. 4 of \cite{Y}. The reason is that Yang \cite{Y}
considers expansions of arbitrary continuous functions on $O$
(not just $\F$-linear ones) with respect to a certain basis $\{
G_n\}$, such that $f_n=G_{q^n}$. Therefore $q^{-n}$ appears in
(8), instead of $n^{-1}$ in the formula from \cite{Y}.

Since the operator $\Delta =\tau d$ is linear,
we have $\Delta f_i=D_i^{-1}\Delta e_i$, $i\ge 1$; clearly
$\Delta f_0=0$. It is known \cite{G1} that
$$
\Delta e_i=\frac{D_i}{D_{i-1}^q}e_{i-1}^q,\quad e_{i-
1}^q=e_i+D_{i-1}^{q-1}e_{i-1}.
$$
Since $D_i=[i]D_{i-1}^q$, we find that
\begin{equation}
\Delta f_i=[i]f_i+f_{i-1},\quad i\ge 1.
\end{equation}

It follows from (6) and (9) that
$$
\Delta u(t)=\sum\limits_{j=0}^\infty (c_{j+1}+[j]c_j)f_j(t).
$$
Substituting into (4) and using uniqueness of the Carlitz
expansion we find a recurrence relation
$$
c_{j+1}+[j]c_j=\lambda c_j,\quad j=0,1,2,\ldots ,
$$
whence, given $c_0$, the solution is determined uniquely by
$$
c_n=c_0\prod\limits_{j=0}^{n-1}(\lambda -[j]).
$$

Suppose that $|\lambda |\ge 1$. Since $|[j]|=q^{-1}$ for $j\ge
1$, we see that $|c_n|=|c_0|\cdot |\lambda |^n\nrightarrow 0$ if
$c_0\ne 0$. This contradiction shows that the equation (4) has no
continuous solutions if $|\lambda |\ge 1$. Therefore we shall
assume that $|\lambda |<1$. Let $u(t,\lambda )$ be the solution
of (4) with $c_0=1$; note that the fixation of $c_0$ is
equivalent to the initial condition $u(1,\lambda )=c_0$. The
function $u(t,\lambda )$ is a function field counterpart of the
power function $t^\lambda$.

\begin{teo}
The function $t\mapsto u(t,\lambda )$, $|\lambda |<1$,
is continuous on $O$. It is
analytic on $O$ if and only if $\lambda =[j]$ for some $j\ge 0$;
in this case $u(t,\lambda )=u(t,[j])=t^{q^j}$. If $\lambda \ne
[j]$ for any integer $j\ge 0$, then $u(t,\lambda )$ is locally
analytic on $O$ if and only if $\lambda =-x$, and in that case
$u(t,-x)=0$ for $|t|\le q^{-1}$. The relation
\begin{equation}
u(t^{q^m},\lambda )=u(t,\lambda^{q^m}+[m]),\quad t\in O,
\end{equation}
holds for all $\lambda$, $|\lambda |<1$, and for all
$m=0,1,2,\ldots$.
\end{teo}

\medskip
{\it Proof}. If $u(t)=t^{q^j}$, $j\ge 0$, then
$$
\Delta u(t)=(xt)^{q^j}-xt^{q^j}=\left( x^{q^j}-x\right)
t^{q^j}=[j]u(t),
$$
so that $u(t,[j])=t^{q^j}$.

Suppose that $\lambda \ne [j]$, $j=0,1,2,\ldots$. Then $|c_n|\le
\{\max (|\lambda |,q^{-1})\}^n\to 0$ as $n\to \infty$, so that
$u(t,\lambda )$ is continuous. More precisely, if $\lambda \ne -x$, then
$|\lambda +x|=q^{-\nu }$ for some $\nu >0$,
$$
|\lambda -[j]|=\left| (\lambda +x)-x^{q^j}\right|=q^{-\nu },\quad j\ge j_0,
$$
if $j_0$ is large enough. This means that for some positive constant $C$
\begin{equation}
|c_n|=Cq^{-n\nu },\quad n\ge j_0.
\end{equation}

On the other hand, if $\lambda =-x$, then $|\lambda -[j]|=q^{-q^j}$,
\begin{equation}
|c_n|=q^{-\frac{q^n-1}{q-1}}.
\end{equation}

If $\lambda \ne -x$, then by (11) $\gamma =0$, so that $u(t,\lambda )$
is not locally analytic. If $\lambda =-x$, we see from (8) and
(12) that $\gamma =(q-1)^{-1}$, $l=1$, and $u(t,-x)$ is analytic
on any ball of the radius $q^{-1}$. We have
$$
u(t,-x)=\sum\limits_{n=0}^\infty (-1)^nx^{\frac{q^n-1}{q-
1}}f_n(t),
$$
and $u(t,-x)$ is not the identical zero on $O$ due to the
uniqueness of the Fourier-Carlitz expansion.

At the same time, since $u(t,-x)$ is analytic on the ball $\{
|t|\le q^{-1}\}$, we can write
$$
u(t,-x)=\sum\limits_{m=0}^\infty a_mt^{q^m},\quad |t|\le q^{-1}.
$$
Substituting this into the equation (4) with $\lambda =-x$, we
find that $a_m=0$ for all $m$, that is $u(t,-x)=0$ for $|t|\le
q^{-1}$.

In order to prove (10), note first that (10) holds for $\lambda
=[j]$, $j=0,1,2,\ldots$. Indeed,
$$
u(t^{q^m},[j])=\left( t^{q^m}\right)^{q^j}=t^{q^{m+j}}
$$
and
$$
[j]^{q^m}+[m]=\left( x^{q^j}-x\right)^{q^m}+x^{q^m}-x=[m+j].
$$

Let us fix $t\in O$. We have
$$
u(t,\lambda )=\sum\limits_{n=0}^\infty \left\{
\prod\limits_{j=0}^{n-1}(\lambda -[j])\right\} f_n(t),
$$
and the series converges uniformly with respect to $\lambda \in
\overline{P}_r$ where
$$
\overline{P}_r=\left\{ \lambda \in \K:\ |\lambda \le r\right\},
$$
for any positive $r<1$. Thus $u(t,\lambda )$ is an analytic
element on $\overline{P}_r$ (see Chapter 10 of \cite{E}). Similarly,
$u(t^{q^m},\lambda )$ and $u(t,\lambda^{q^m}+[m])$ are analytic
elements on $\overline{P}_r$ (for the latter see Theorem 11.2
from \cite{E}). Suppose that $q^{-1}\le r<1$. Since both sides of
(10) coincide on an infinite sequence of points $\lambda =[j]$,
$j=0,1,2,\ldots$, they coincide on $\overline{P}_r$ (see
Corollary 23.8 in \cite{E}). This implies their coincidence for
$|\lambda |<1$. $\quad \blacksquare$

\bigskip
Similarly, if in (4) $\lambda$ is a $m\times m$ matrix with
elements from $\K$ (we shall write $\lambda \in M_m(\K)$), and we
look for a solution $u\in M_m(\K)$, then we can find a continuous
solution (6) with matrix coefficients
\begin{equation}
c_i=\left\{ \prod\limits_{j=0}^{i-1}(\lambda -[j]I_m)\right\}
c_0,\quad i\ge 1
\end{equation}
($I_m$ is a unit matrix) if $|\lambda
|\stackrel{\mbox{{\footnotesize def}}}{=}\max |\lambda_{ij}|<1$. Note that
$c_0=u(1)$, so that if $c_0$ is an invertible matrix, then $u$ is
invertible on a certain neighbourhood of 1.

\section{FIRST ORDER SYSTEMS}

\medskip
Let us consider a system (5) with the coefficient $P(\tau )$
given in (2). We assume that $|\pi_k|\le \gamma$, $\gamma >0$,
for all $k$, $|\pi_0|<1$. Denote by $g(t)$ a solution of the
equation $\tau dg=\pi_0g$. Let $\lambda_1,\ldots ,\lambda_m\in
\K$ be the eigenvalues of the matrix $\pi_0$.

\begin{teo}
If
\begin{equation}
\lambda_i-\lambda_j^{q^k}\ne [k],\quad i,j=1,\ldots ,m;\
k=1,2,\ldots,
\end{equation}
then the system (5) has a matrix solution
\begin{equation}
u(t)=W(g(t)),\quad W(s)=\sum\limits_{k=0}^\infty w_ks^{q^k},\quad
w_0=I_m,
\end{equation}
where the series for $W$ has a positive radius of convergence.
\end{teo}

{\it Proof}. Substituting (15) into (5), using the fact that
$\Delta =\tau d$ is a derivation of the composition ring of
$\F$-linear series, and that $\Delta (t^{q^k})=[k]t^{q^k}$, we
come to the identity
\begin{equation}
\sum\limits_{k=0}^\infty [k]w_k\tau^k(g(t))+\sum\limits_{k=0}^\infty
w_k\tau^k(\pi_0)\tau^k(g(t))
-\sum\limits_{j=0}^\infty \pi_j\tau^j\left(
\sum\limits_{k=0}^\infty w_k\tau^k(g(t))\right) =0.
\end{equation}

If the series for $W$ has indeed a positive radius of convergence
(which will be proved later), then all the expressions in (16)
make sense for a small $|t|$, since $g(t)\to 0$ as $|t|\to 0$.
Since $w_0=I_m$, the first summand in the second sum in (16) and
the summand with $j=k=0$ in the third sum are cancelled. Changing
the order of summation we find that (16) is equivalent to the
system of equations
\begin{equation}
w_k\left( [k]I_m+\tau^k(\pi_0)\right) -
\pi_0w_k=\sum\limits_{j=1}^k\pi_j\tau^j(w_{k-j}),\quad
k=1,2,\ldots ,
\end{equation}
with respect to the matrices $w_k$.

The system (17) is solved step by step -- if the right-hand side
of (17) with some $k$ is already known, then $w_k$ is determined
uniquely, provided the spectra of the matrices
$[k]I_m+\tau^k(\pi_0)$ and $\pi_0$ have an empty intersection
(\cite{Ga}, Sect. VIII.1). This condition is equivalent to
(14), and it remains
to prove that the series for $W$ has a non-zero radius of
convergence.

Let us transform $\pi_0$ to its Jordan normal form. We have $U^{-
1}\pi_0U=A$ where $U$ is an invertible matrix over $\K$, and $A$
is block-diagonal:
$$
A=\bigoplus\limits_{\alpha =1}^l\left( \lambda^{(\alpha
)}I_{d_\alpha }+H^{(\alpha )}\right)
$$
where $\lambda^{(\alpha )}$ are eigenvalues from the collection
$\{ \lambda_1,\ldots ,\lambda_m\}$, $H^{(\alpha )}$ is a Jordan
cell of the order $d_\alpha$ having zeroes on the principal
diagonal and 1's on the one below it. Denote $\mu_k^{(\alpha
)}=\left( \lambda^{(\alpha )}\right)^{q^k}+[k]$. If
$V_k=\tau^k(U)$, then
\begin{equation}
V_k^{-1}\left( [k]I_m+\tau^k(\pi_0)\right) V_k=\bigoplus\limits_{\alpha
=1}^l\left( \mu_k^{(\alpha )}I_{d_\alpha }+H^{(\alpha )}\right) .
\end{equation}

If $B_k$ is the matrix (18), and $C_k$ is the matrix in the
right-hand side of (17), then (17) takes the form
$$
w_kV_kB_kV_k^{-1}-UAU^{-1}w_k=C_k
$$
or, if we use the notation $\widetilde{w}_k=U^{-1}w_kV_k$,
\begin{equation}
\widetilde{w}_kB_k-A\widetilde{w}_k=\widetilde{C}_k,\quad
k=1,2,\ldots,
\end{equation}
where
$$
\widetilde{C}_k=U^{-1}C_kV_k=U^{-1}\left(
\sum\limits_{j=1}^k\pi_j\tau^j\left( U\widetilde{w}_{k-j}V_{k-
j}^{-1}\right) \right) V_k=U^{-
1}\sum\limits_{j=1}^k\pi_j\tau^j(U)\tau^j\left( \widetilde{w}_{k-
j}\right) ,
$$
$\widetilde{w}_0=I_m$. We may assume that $|U|\le 1$, $|U^{-
1}|\le \rho$, $\rho >0$.

In accordance with the quasi-diagonal form of the matrices $A$
and $B_k$ we can decompose the matrix $\widetilde{w}_{k}$ into
$d_\alpha \times d_\beta$ blocks
$$
\widetilde{w}_{k}=\left( \widetilde{w}_{k}^{(\alpha \beta
)}\right) ,\quad \alpha ,\beta =1,\ldots ,l.
$$
Similarly we write $\widetilde{C}_{k}=\left( \widetilde{C}_{k}^{(\alpha
\beta )}\right)$. Then the system (19)
is decoupled into a system of equations for each block:
\begin{equation}
\left( \mu_k^{(\beta )}-\lambda^{(\alpha )}\right)
\widetilde{w}_{k}^{(\alpha \beta )}-H^{(\alpha )}
\widetilde{w}_{k}^{(\alpha \beta )}+\widetilde{w}_{k}^{(\alpha \beta
)}H^{(\beta )}=\widetilde{C}_{k}^{(\alpha \beta )}.
\end{equation}

The equation (20) can be considered as a system of scalar
equations with respect to elements of the matrix
$\widetilde{w}_{k}^{(\alpha \beta )}$. Let us enumerate these
elements $\left( \widetilde{w}_{k}^{(\alpha \beta )}\right)_{ij}$
lexicographically (in i,j) with the inverse enumeration order of
the second index $j$. The product $H^{(\alpha )}
\widetilde{w}_{k}^{(\alpha \beta )}$ is obtained from
$\widetilde{w}_{k}^{(\alpha \beta )}$ by the shift of all the
rows one step upwards, the last row being filled by zeroes.
Similarly, the product $\widetilde{w}_{k}^{(\alpha \beta
)}H^{(\beta )}$ is the result of shifting all the columns of
$\widetilde{w}_{k}^{(\alpha \beta )}$ one step to the right and
filling the first column by zeroes (\cite{Ga}, Sect. I.3).
This means that the
system (2) (with fixed $\alpha ,\beta$) with the above
enumeration of the unknowns is upper triangular. Indeed, the
latter is equivalent to the fact that each equation contains,
together with some unknown, only the unknowns with larger
numbers, and this property is the result of the above shifts.

Therefore the determinant $D_k^{(\alpha \beta)}$ of the system
(20) equals $\left( \mu_k^{(\beta )}-\lambda^{(\alpha
)}\right)^{d_\alpha d_\beta }$. By our assumption $|\pi_0|<1$,
and if $\lambda^{(\alpha )}$ is an eigenvalue of $\pi_0$ with an
eigenvector $f\ne 0$, then $|\lambda^{(\alpha )}|\cdot
|f|=|\pi_0f|<|f|$, so that $|\lambda^{(\alpha )}|<1$. This means
that all the coefficients on the left in (2) have the absolute
values $\le 1$.

It follows from (14) that $\lambda_i\ne -x$, $i=1,\ldots ,n$. As
$k\to \infty$, $\mu_k^{(\beta )}=\left( \lambda^{(\beta
)}\right)^{q^k}+x^{q^k}-x\to -x$. Thus $\left|
\mu_k^{(\beta )}-\lambda^{(\alpha )}\right| \ge \sigma_1>0$ for
all $k$, whence $\left| D_k^{(\alpha \beta)}\right| \ge \sigma_2>0$
where $\sigma_2$ does not depend on $k$. Now we obtain an
estimate for the solution of the system (19),
\begin{equation}
\left| \widetilde{w}_{k}\right| \le \rho_1\left|
\widetilde{C}_{k}\right|,
\end{equation}
with $\rho_1>0$ independent of $k$.

Looking at (21) we find that
$$
\left| \widetilde{w}_{k}\right| \le \rho_2\max\limits_{1\le j\le
k}\left| \widetilde{w}_{k-j}\right|^{q^j}
$$
where $\rho_2$ does not depend on $k$. We may assume that
$\rho_2\ge 1$. Now we find that
\begin{equation}
\left| \widetilde{w}_{k}\right| \le \rho_2^{q^{k-1}+q^{k-
2}+\cdots +1},\quad k=1,2,\ldots .
\end{equation}

Indeed, (22) is obvious for $k=1$. Suppose that we have proved
the inequalities
$$
\left| \widetilde{w}_{j}\right| \le \rho_2^{q^{j-1}+q^{j-
2}+\cdots +1},\quad 1\le j\le k-1.
$$
Then
\begin{multline*}
\left| \widetilde{w}_{k}\right| \le \rho_2\max \left( 1,\left|
\widetilde{w}_{1}\right|^{q^{k-1}},\ldots ,\left|
\widetilde{w}_{k-1}\right|^q\right) \\
\le \rho_2\max \left( 1,\rho_2^{q^{k-1}},\rho_2^{(q+1)q^{k-2}},
\ldots ,\rho_2^{(q^{k-2}+\cdots +1)q}\right) =\rho_2^{q^{k-
1}+q^{k-2}+\cdots +1},
\end{multline*}
and (22) is proved.

Therefore, since $w_k=U\widetilde{w}_k\tau^k(U^{-1})$, we have
$$
|w_k|\le \rho^{q^k}\cdot \rho_2^{q^{k-1}+\cdots +1}\le
\rho_3^{\frac{q^{k+1}-1}{q-1}},\quad \rho_3>0,
$$
which means that the series in (15) has a positive radius of
convergence. $\quad \blacksquare$

\bigskip
{\it Remarks}. 1). If $\varphi \in \mathbf F_q^m$, then, as usual,
$v=u\varphi$ is a vector solution of the system $\tau dv-P(\tau
)v=0$, since the system is $\F$-linear. However, the system is
nonlinear over $\K$, so that we cannot obtain a vector solution
in such a way for an arbitrary $\varphi \in \K$.

2). Analogs of the condition (14) occur also in the analytic
theory of differential equations over $\mathbb C$ (see Corollary
11.2 in \cite{H}) and $\mathbb Q_p$ (Sect. III.8 in \cite{DGS}).
For systems over $\mathbb C$, it is requested that differences of
the eigenvalues of the leading coefficient $\pi_0$ must not be
non-zero integers. Over $\mathbb Q_p$, in addition to that,
the eigenvalues must not be non-zero integers themselves.
In both the characteristic zero cases it is possible to get rid
of such conditions by using special changes of variables called
shearing transformations. For example, let $m=1$, and the equation over
$\mathbb Q_p$ has the form
$$
\zeta u'(\zeta )=nu(\zeta )+\left( \sum\limits_{k=1}^\infty
\pi_k\zeta^k\right) u(\zeta ),\quad n\in \mathbb N.
$$
Then the change of variables $u(\zeta )=\zeta v(\zeta )$ gives the
transformed equation
$$
\zeta v'(\zeta )=(n-1)v(\zeta )+\left( \sum\limits_{k=1}^\infty
\pi_k\zeta^k\right) v(\zeta ).
$$
Repeating the transformation, we remove the term violating the
condition. A modification of this approach works for systems of
equations.

In our case the situation is different.
Let us consider again the scalar case $m=1$. Here the condition (14) is
equivalent to the condition $\pi_0\ne -x$ (the general solution
of the equation $\pi_0-\pi_0^q=[k]$ has the form $\pi_0=-x+\xi$ where
$\xi -\xi^{q^k}=0$, that is either $\xi =0$, or $|\xi |=1$; the latter
contradicts our assumption $|\pi_0|<1$). If, on the contrary, $\pi_0=-x$,
then, as we saw in Theorem 1, $g(t)=0$ for $|t|\le q^{-1}$, and the
construction (15) does not make sense. On the other hand, a formal
analog of the shearing transformation for this case is the
substitution $u=\tau (v)$. However it is easy to see
that $v$ satisfies an equation with the same principal part,
as the equation for $u$.

\section{EULER TYPE EQUATIONS}

\medskip
Classically, the Euler equation has the form
$$
\zeta^mu^{(m)}(\zeta )+\beta_{m-1}\zeta^{m-1}u^{(m-1)}(\zeta
)+\cdots +\beta_0u(\zeta )=0
$$
where $\beta_0,\beta_1,\ldots \beta_{m-1}\in \mathbb C$. It can
be reduced to a first order linear system with a constant matrix.
The solutions are linear combinations of functions of the form
$\zeta^\lambda (\log \zeta )^k$. Of course, such functions have
no direct $\F$-linear counterparts, and our study of the Euler-type
equations will again be based on expansions in the Carlitz
polynomials.

Let us consider a linear equation
\begin{equation}
\tau^md^mu+b_{m-1}\tau^{m-1}d^{m-1}u+\cdots +b_0u=0
\end{equation}
where $b_0,b_1,\ldots ,b_{m-1}\in \K$. In order to reduce (23) to
a first order system, it is convenient to set
$$
\varphi_k=\tau^{k-1}d^{k-1}u,\quad k=1,\ldots ,m.
$$
Since $d\tau^{k-1}-\tau^{k-1}d=[k-1]^{1/q}\tau^{k-2}$ (see
\cite{K3}), we have
$$
\tau d\varphi_k=\tau^kd^ku+[k-1]\tau^{k-1}d^{k-
1}u=\varphi_{k+1}+[k-1]\varphi_k,
$$
$k=1,\ldots ,m-1$. Next, by (23),
$$
\tau d\varphi_m=\tau^md^mu+[m-1]\tau^{m-1}d^{m-1}u=([m-1]-b_{m-
1})\varphi_m-b_{m-2}\varphi_{m-1}-\cdots -b_0\varphi_1.
$$

Thus the equation (23) can be written as a system
\begin{equation}
\tau d\varphi =B\varphi ,\quad \varphi =(\varphi_1,\ldots
,\varphi_m),
\end{equation}
where
$$
B=
\begin{pmatrix}0 & 1 & 0 & 0 & \ldots & 0\\
0 & [1] & 1 & 0 & \ldots & 0\\
0 & 0 & [2] & 1 & \ldots & 0\\
0 & 0 & 0 & [3] & \ldots & 0\\
\hdotsfor[1]{6}\\
0 & 0 & 0 & 0 & \ldots & 1\\
-b_0 & -b_1 & -b_2 & -b_3 & \ldots & [m-1]-b_{m-1}\end{pmatrix}.
$$
This time we cannot use directly the above results, since $|B|\ge
1$. However in some cases it is possible to proceed in a slightly
different way.

Suppose that all the eigenvalues of the matrix $B$ lie in the
open disk $\{ |\lambda |<1\}$. Transforming $B$ to its Jordan
normal form we find that
$$
B=X^{-1}(B_0+N)X
$$
where $X$ is an invertible matrix, $B_0$ is a diagonal matrix,
$|B_0|=\mu <1$, $N$ is nilpotent, that is $N^\varkappa =0$ for
some natural number $\varkappa$, and $N$ commutes with $B_0$. If
$\Psi$ is a matrix solution of the system
$$
\tau d\Psi =(B_0+N)\Psi,
$$
then $\Phi =X^{-1}\Psi X$ is a matrix solution of (24).

On the other hand, we can obtain $\Psi$ just as in the case $N=0$
considered in Sect. 2, writing
\begin{equation}
\Psi (t)=\sum\limits_{i=0}^\infty c_if_i(t),
\end{equation}
\begin{equation}
c_i=\left\{ \prod\limits_{j=0}^{i-1}(B_0+N-[j]I_m)\right\}c_0.
\end{equation}

Indeed, the product in (26) is the sum of the expressions $(-
[j])^{\nu_1}B_0^{\nu_2}N^{\nu_3}$ where $\nu_1+\nu_2+\nu_3=i$,
$\nu_3<\varkappa$. Therefore in (25)
$$
|c_i|\le |c_0|\cdot |N|^{\varkappa -1}\left\{ \max (\mu ,
q^{-1})\right\}^{i-\varkappa }\longrightarrow 0,\quad i\to \infty .
$$

Let us consider in detail the case where $m=2$. Our equation is
\begin{equation}
\tau^2d^2u+b_1\tau du+b_0u=0.
\end{equation}
Now we have the system (24) with
$$
B=\begin{pmatrix}
0 & 1\\
-b_0 & [1]-b_1\end{pmatrix}.
$$

The characteristic polynomial of $B$ is $D_2(\lambda )=\lambda^2+\lambda
(b_1-[1])+b_0$, with the discriminant $\delta =(b_1-[1])^2-4b_0$. We
assume that the eigenvalues are such that
$|\lambda_1|,|\lambda_2|<1$. This condition is satisfied, for
example, if $p\ne 2$, $|b_0|<1$, $|b_1|<1$.

The greatest common divisor of the first order minors of $B-
\lambda I_2$ is 1. This means that $B$ is diagonalizable if and
only if $\lambda_1\ne \lambda_2$, that is if $\delta \ne 0$ (see
e.g. \cite{Gel}). In this case
\begin{equation}
B=X^{-1}\begin{pmatrix}
\lambda_1 & 0\\
0 & \lambda_2\end{pmatrix}X
\end{equation}
for some invertible matrix $X$, and our system has a matrix
solution $\Phi$, such that
$$
X\Phi (t)X^{-1}=tI_2+\sum\limits_{n=1}^\infty \diag \left\{
\prod\limits_{j=0}^{n-1}(\lambda_1-[j]),\prod\limits_{j=0}^{n-1}
(\lambda_2-[j])\right\} f_n(t)\stackrel{\mbox{{\footnotesize def}}}{=}\diag
\{
\psi_1(t),\psi_2(t)\}.
$$

It is easy to see that $\psi_1(t)$ and $\psi_2(t)$ are solutions
of the equation (27). Indeed, if $X^{-1}=\left( \begin{smallmatrix}
\xi_{11} & \xi_{12}\\ \xi_{21} & \xi_{22}\end{smallmatrix}\right)$, then
$$
\Phi (t)X^{-1}\begin{pmatrix}1\\0\end{pmatrix}=(\xi_{11}\psi_1(t),
\xi_{21}\psi_2(t)),
$$
and (for $\psi_1$) it is sufficient to show that $\xi_{11}\ne 0$.
However $X^{-1}X=I_2$, and if $\xi_{11}=0$, then writing
$X=\left( \begin{smallmatrix}
\chi_{11} & \chi_{12}\\ \chi_{21} & \chi_{22}\end{smallmatrix}\right)$
we find that $\xi_{12}\chi_{22}=0$, $\xi_{12}\chi_{21}=1$. At the
same time, by (28), $\xi_{12}\lambda_2\chi_{21}=0$, and
$\xi_{12}\lambda_2\chi_{22}=1$, and we come to a contradiction. A
similar reasoning works for $\psi_2(t)$. It follows from the
uniqueness of the Fourier-Carlitz expansion that $\psi_1$ and $\psi_2$ are
linearly independent.

If $\lambda_1=\lambda_2\stackrel{\mbox{{\footnotesize def}}}{=}\lambda$,
then $B$ is similar to the Jordan cell
$$
N=\begin{pmatrix}
\lambda & 1\\
0 & \lambda \end{pmatrix}.
$$
It is proved by induction that
$$
\prod\limits_{j=0}^{n-1}(N-[j]I_2)=\begin{pmatrix}
\prod\limits_{j=0}^{n-1}(\lambda -[j]) & \sum\limits_{j=0}^{n-1}
\prod\limits_{\substack{0\le i\le n-1\\i\ne j}}(\lambda -[i])\\
0 & \prod\limits_{j=0}^{n-1}(\lambda -[j])\end{pmatrix}.
$$
In this case we have the following two linearly independent solutions of
(27):
$$
\psi_1(t)=t+\sum\limits_{n=1}^\infty
\left\{ \prod\limits_{j=0}^{n-1}(\lambda -[j])\right\} f_n(t),
$$
$$
\psi_2(t)=t+\sum\limits_{n=1}^\infty \left\{ \sum\limits_{j=0}^{n-1}
\prod\limits_{\substack{0\le i\le n-1\\i\ne j}}(\lambda -[i])
\right\} f_n(t).
$$

Thus, for the case of the eigenvalues from the disk $\{ |\lambda |<1\}$, we
have given an explicit construction of solutions for the Euler
type equations.

\section{DISCONTINUOUS SOLUTIONS}
\medskip
For all the above equations, the solutions were found as Fourier-Carlitz
expansions (6), and we had to impose certain conditions upon
coefficients of the equations, in order to guarantee the uniform
convergence of the series (6) on $O$. However, formally we could
write the series for the solutions without those conditions. Thus
it is natural to ask whether the corresponding series (6)
converge at some points $t\in O$. Note that (6) always makes
sense for $t\in \F[x]$ (for each such $t$ only a finite number of
terms is different from zero). The question is whether the series
(6) converges on a wider set; if the answer is negative, such a
formal solution will be called {\it strongly singular}.

We shall need the following property of the Carlitz polynomials.

\begin{lem}
$$
\left| f_i(x^n)\right| =\begin{cases}
0, & \text{if $n<i$},\\
q^{i-n}, & \text{if $n\ge i$}.\end{cases}
$$
\end{lem}

{\it Proof}. If $n<i$, it follows from (7) that $f_i(x^n)=0$. Let
$n\ge i$. Then $|x^n-\omega |=|\omega |$ for all $\omega \in \F
[x]$, $\deg \omega <i$. Writing
$$
e_i(t)=t\prod\limits_{\substack{0\ne w\in \F [x]\\ \deg \omega
<i}}(t-\omega )
$$
we find that
$$
\left| e_i(x^n)\right| =|x^n|\prod\limits_{\substack{\deg \omega
<i\\ \omega \ne 0}}|\omega |=|x^n|\cdot \left| \lim\limits_{t\to
0}\frac{e_i(t)}{t}\right| .
$$

It is known \cite{C1,G1,G2} that
$$
e_i(t)=\sum\limits_{j=0}^i(-1)^{i-j}\frac{D_i}{D_jL_{i-
j}^{q^i}}t^{q^j}
$$
where $L_i=[i][i-1]\cdots [1]$, $i\ge 1$, and $L_0=1$. Therefore
$$
\lim\limits_{t\to 0}\frac{e_i(t)}{t}=(-1)^i\frac{D_i}{L_i}
$$
whence
$$
\left| f_i(x^n)\right| =\frac{q^{-n}}{|L_i|}=q^{i-n}
$$
as desired. $\quad \blacksquare$

\medskip
Now we get a general sufficient condition for a function (6) to
be strongly singular.

\begin{teo}
If $|c_i|\ge \rho >0$ for all $i\ge i_0$ (where $i_0$ is some
natural number), then the function (6) is strongly singular.
\end{teo}

{\it Proof}. In view of the convergence criterion for series in a
non-Archimedean field (see Sect. 1.1.8 in \cite{BGR}), it is
sufficient to find, for any $t\ne \F [x]$, such a sequence
$i_k\to \infty$ that $\left| f_{i_k}(t)\right| =1$, $k=1,2,\ldots
$.

In fact, if $t\in O\setminus \F [x]$, then
$t=\sum\limits_{n=0}^\infty \xi_nx^n$, $\xi_n\in \F$, with
$\xi_{i_k}\ne 0$ for some sequence
$i_k\to \infty$. We have (by Lemma 1)
$$
f_{i_k}(t)=\sum\limits_{n=i_k}^\infty \xi_nf_{i_k}(x^n)
$$
where $\left| f_{i_k}(x^{i_k})\right|=1$, $\left| \xi_{i_k}\right| =1$,
$\left| f_{i_k}(x^n)\right| =q^{i_k-n}\le q^{-1}$ for $n>i_k$.
Thus $\left| f_{i_k}(t)\right| =1$ for all $k$. $\quad
\blacksquare$

\medskip
It follows from Theorem 3 and the discussion preceding Theorem 1
that non-trivial formal solutions of
the equation (4) with $|\lambda |\ge 1$ are strongly singular. A
more complicated example of an equation with such solutions will
be given in the next section.

\section{HYPERGEOMETRIC EQUATION}

\medskip
The equation for Thakur's function field analog of the
hypergeometric function ${}_2F_1(a,b;1;t)$ has the form
\begin{equation}
(\Delta -[-a])(\Delta -[-b])u=d\Delta u
\end{equation}
where $a,b\in \mathbb Z$.

A corresponding classical equation over $\mathbb C$ has the form
$$
\left( \zeta \frac{d}{d\zeta}+a\right) \left( \zeta
\frac{d}{d\zeta}+b\right) u=\frac{d}{d\zeta}\left( \zeta
\frac{du}{d\zeta}\right) .
$$
Its holomorphic solution near the origin is ${}_2F_1(a,b;1;\zeta
)$; the second solution has a logarithmic singularity (see
\cite{Olv}, Chapter 5, $\S \S$10, 11).

As before, we look for a solution of the form (6) defined at
least for $t\in \F [x]$; note that the operators $\Delta$ and $d$
are well-defined on functions of $t\in \F [x]$. Using the
relation (8) and the fact that $df_0=0$, $df_i=f_{i-1}$ ($i\ge
1$) we obtain a recursive relation
\begin{multline}
\left( c_{i+2}^{1/q}-c_{i+2}\right) +c_{i+1}^{1/q}[i+1]^{1/q}-
c_{i+1}([i]+[i+1]-[-a]-[-b])\\
-c_i([i]-[-a])([i]-[-b])=0,\quad i=0,1,2,\ldots .
\end{multline}

Taking arbitrary initial coefficients $c_0,c_1\in \K$ we obtain a
solution $u$ defined on $\F [x]$. On each step we have to solve
the equation
\begin{equation}
z^{1/q}-z=v.
\end{equation}
If $|c_i|\le 1$ and $|c_{i+1}|\le 1$, then in the equation for
$c_{i+2}$ we have $|v|<1$.

\begin{lem}
The equation (31) with $|v|<1$ has a unique solution $z_0\in \K$,
for which $|z_0|\le |v|$, and $q-1$ other solutions $z$, $|z|=1$.
\end{lem}

{\it Proof}. It is convenient to investigate the equivalent
equation
$$
z^q-z=w,\quad |w|<1.
$$

Consider the polynomial $\varphi (z)=z^q-z-w$, $\varphi \in
O[z]$. Since $\varphi'(z)\equiv -1$ and $\varphi (w^{1/q})=-
w^{1/q}$, we have $|\varphi (w^{1/q})|=|w|^{1/q}<1=\left| \left\{
\varphi'(w^{1/q})\right\}^2\right|$. Thus we are within the
conditions of the version of the Hensel lemma for a field with a
possibly non-discrete valuation, given in Chapter XII of
\cite{L}. By that result the polynomial $\varphi$ has such a root
$z_0\in \K$ that $\left| z_0-w^{1/q}\right| \le |w|^{1/q}$,
whence
$$
|z_0|\le \max \left( \left| z_0-w^{1/q}\right|,|w|^{1/q}\right) =
|w|^{1/q}.
$$
Taking $w=-v^q$ we obtain the required solution of (31).

Other solutions of (31) have the form $z=z_0+\theta$, $0\ne
\theta \in \F$. Obviously $|z|=1.\quad \blacksquare$

\medskip
It follows from Lemma 2 that if $|c_i|,|c_{i+1}|\le 1$ for some
$i$, then $|c_n|\le 1$ for all $n\ge i$.

The relation (30), Lemma 2, and Theorem 3 imply the following
property of solutions of the equation (29). It is natural to call
a solution {\it generic} if, starting from a certain step of
finding the coefficients $c_n$, we always take the most frequent
option corresponding to a solution of (31) with $|z|=1$.

\begin{teo}
A generic solution of the equation (29) is strongly singular.
\end{teo}

\medskip
Of course, in some special cases the recursion (30) can lead to
more regular solutions, in particular, to the holomorphic
solutions found by Thakur \cite{Th1,Th2}.


\medskip
\begin{thebibliography}{999}
\bibitem{BGR}
S. Bosch, U. G\"untzer, and R. Remmert, ``Non-Archimedean
Analysis'', Springer, Berlin, 1984.
\bibitem{C1}
L. Carlitz, On certain functions connected with polynomials in a
Galois field, {\it Duke Math. J.} {\bf 1} (1935), 137--168.
\bibitem{C2}
L. Carlitz, Some special functions over $GF(q,x)$, {\it Duke Math. J.}
{\bf 27} (1960), 139--158.
\bibitem{CL}
E. A. Coddington and N. Levinson, ``Theory of Ordinary
Differential Equations'', McGraw-Hill, New York, 1955.
\bibitem{DGS}
B. Dwork, G. Gerotto, and F. J. Sullivan, ``An Introduction to
$G$-Functions'', Princeton University Press, 1994.
\bibitem{E}
A. Escassut, ``Analytic Elements in $p$-Adic Analysis'', World
Scientific, Singapore, 1995.
\bibitem{Ga}
F. R. Gantmacher, ``Matrizentheorie'', Springer, Berlin, 1986.
\bibitem{Gel}
I. M. Gelfand, ``Lectures on Linear Algebra'', Interscience, New
York, 1961.
\bibitem{G1}
D. Goss, Fourier series, measures, and divided power series in
the theory of function fields, {\it K-Theory} {\bf 1} (1989),
533--555.
\bibitem{G2}
D. Goss, ``Basic Structures of Function Field Arithmetic'',
Springer, Berlin, 1996.
\bibitem{H}
P. Hartman, ``Ordinary Differential Equations'', Wiley, New York,
1964.
\bibitem{K1}
A. N. Kochubei, Harmonic oscillator in characteristic $p$, {\it Lett.
Math. Phys.} {\bf 45} (1998), 11--20.
\bibitem{K2}
A. N. Kochubei, $\F$-linear calculus over function fields, {\it J.
Number Theory} {\bf 76} (1999), 281--300.
\bibitem{K3}
A. N. Kochubei, Differential equations for $\F$-linear functions, {\it J.
Number Theory} {\bf 83} (2000), 137--154.
\bibitem{L}
S. Lang, ``Algebra'', Addison-Wesley, Reading, 1965.
\bibitem{Olv}
F. W. J. Olver, ``Asymptotics and Special Functions'', Academic
Press, New York, 1974.
\bibitem{RC}
P. Robba and G. Christol, ``\'Equations Diff\'erentielles
$p$-Adiques'', Hermann, Paris, 1994.
\bibitem{SL}
D. Sinnou and D. Laurent, Ind\'ependence algebrique sur les
$T$-modules, {\it Compositio Math.} {\bf 122} (2000), 1--22.
\bibitem{Th1}
D. Thakur, Hypergeometric functions for function fields, {\it Finite
Fields and Their Appl.} {\bf 1} (1995), 219--231.
\bibitem{Th2}
D. Thakur, Hypergeometric functions for function fields II,
{\it J. Ramanujan Math. Soc.} {\bf 15} (2000), 43--52.
\bibitem{Y}
Z. Yang, Locally analytic functions over completions of $\mathbf
F_r(U)$, {\it J. Number Theory} {\bf 73} (1998), 451--458.
\end{thebibliography}
\end{document}